\documentclass[preprint,12pt]{elsarticle}

\usepackage{amssymb}

\journal{Examples and Counterexamples}
\usepackage{xcolor}
\definecolor{darkgreen}{rgb}{0.0,0.6,0.0}
\newtheorem{theorem}{Theorem}
\newtheorem{lemma}{Lemma}
\newtheorem{definition}{Definition}

\begin{document}

\begin{frontmatter}



\title{Infinitesimal phase response functions can be misleading}


\author{Christoph B\"orgers}

\affiliation{organization={Department of Mathematics, Tufts University}, \\
            city={Medford, MA 02421,}, 
            country={USA}}

\begin{abstract}
Phase response functions are the central tool in the mathematical analysis of pulse-coupled oscillators. When an 
oscillator receives a brief input pulse, the phase response function 
specifies how its phase shifts as a function of  the phase at which the input is received. When the pulse is weak, it is customary 
to linearize around zero pulse strength.
The result is called the {\em infinitesimal} phase response function. 
These ideas have been used extensively
in theoretical biology, and also in some areas of engineering. I give  examples showing that 
the infinitesimal phase response function may predict that two oscillators, as they exchange pulses back and fourth, will 
converge to synchrony, yet this is false when the exact phase response function is used, for all positive
interaction strengths.
For short, the analogue of the 
Hartman-Grobman theorem that one might expect to hold at first sight is invalid. 
I give a condition under which the prediction derived using the infinitesimal phase
response function does hold for the exact phase response function when interactions
are sufficiently weak but positive. However, I argue 
that this condition may often fail to hold. 
\end{abstract}




\end{frontmatter}


\section{Phase oscillators, phase response functions, and synchrony} {

In this paper, I present a greatly expanded discussion of a calculation  from \cite{borgers_book}, concerning pulse-coupled phase oscillators. Although I have  in mind periodically firing
neurons, or perhaps populations of neurons that periodically fire spike volleys, the mathematics may be applicable in other 
contexts  as well. In this section, I review some  background.

\subsection{Phase oscillators} 

\vskip 5pt

Stable oscillators are usually described mathematically as points moving along an attracting limit cycle 
in a phase space. For simplicity, we will say that the {\em oscillator} (rather than the point in
phase space) moves along the limit cycle. 

A {\em phase oscillator} is one that never leaves the limit cycle, not even in response to inputs.
Further, the dynamics are assumed to be {\em autonomous}, meaning that in the absence of inputs, the velocity
of motion only depends on position in phase space.
Therefore at 
any time, even when there have been {\em past} inputs, the remaining time to cycle completion in the absence of {\em future }inputs
lies between $0$ and $T$, where $T$ is the period of the oscillator in the absence of {\em any} inputs. 
``Cycle completion" must, of course, be defined in an arbitrary way. We think of neurons firing
periodic trains of action potentials, or networks of neurons firing periodic spike volleys, 
and mean firing by ``cycle completion". 

We define the {\em phase} of a phase oscillator to be 
$$
\varphi = \frac{T-r}{T}
$$
where $r$ is the remaining time to cycle completion in the absence of future inputs.
We have 
$$
\varphi \in [0,1]
$$
always, and  the state of the oscillator is  completely characterized by $\varphi$. 

The phase description is a common idealization in the study of coupled oscillators. It is  an approximation to reality, always somewhat inaccurate, and often 
inappropriate. For example, an input to a neuron can delay the next firing by far more than $T$ \cite[Figure 1]{phoka2010}, and it can even
halt firing altogether \cite{BKG}. Such examples do involve leaving the
limit cycle, and are outside the scope of this article. 

The main example  in this paper is the Ermentrout-Kopell theta neuron model \cite{ermentroutkopell}, which 
explicitly incorporates the assumption that the oscillator never leaves its limit cycle in 
response to external inputs; in fact, the model describes the motion of a point on the unit circle. 

To view an oscillator as a {\em phase} oscillator  is  a more plausible approximation when inputs to the oscillator
are { weaker} or the limit cycle  underlying the oscillation is more strongly attracting. Nonetheless the phase description is commonly used not only  for  {{\em infinitesimally} weak} inputs, but also 
for inputs of positive strength. In neuroscience, inputs of positive (sometimes 
considerable) strength are certainly of interest, especially if we think of an oscillator as representing 
a network of neurons.

\subsection{Phase response functions}

We assume inputs to the oscillators to be 
 {\em instantaneous pulses}. This, too, is a simplification, and it is commonly denoted with the phrase {\em pulse-coupling}. The simplification is valid perhaps for instance as an idealization of
fast AMPA-receptor-mediated excitatory synaptic pulses, but not valid  in many other context. 

When an oscillator at phase $\varphi$ receives an input pulse of strength $\epsilon$, we assume that its phase changes
instantaneously as a result: 
$$
\varphi \mapsto \varphi+g(\varphi, \epsilon).
$$
The function $g$ is called the {\em phase response function.} Phase response functions have been
measured experimentally in many studies in the neuroscience literature; see for instance \cite{phoka2010}.

We will not yet be specific about what we mean by the ``strength" of the input pulse. The parameter 
$\epsilon$ is assumed to be $\geq 0$, with 
\begin{equation} 
\label{eq:nichts} 
g(\varphi,0) = 0 ~~~\mbox{for all $\varphi \in [0,1]$}.
\end{equation}

\noindent
Since phase, by our definition, always lies between $0$ and $1$, we must have 
\begin{equation}
\label{eq:range_of_phi}
- \varphi \leq g(\varphi, \epsilon) \leq 1 - \varphi ~~~~\mbox{for $\varphi \in [0,1]$, $\epsilon \in [0,\infty)$}. 
\end{equation}
We further assume  that receipt of an input cannot {\em instantly} set a neuron to phase $0$ or phase $1$: 
\begin{equation}
\label{eq:range_of_phi_strict}
- \varphi < g(\varphi, \epsilon) < 1 - \varphi ~~~~\mbox{for $\varphi \in (0,1)$, $\epsilon \in [0,\infty)$}.
\end{equation}
``Phase $0$" and ``phase $1$" denote the same thing --- the time of firing. Nonetheless we think of
$g(0, \epsilon)$ and $g(1,\epsilon)$ as   two different quantities: $g(0,\epsilon)$ is the phase shift that an input immediately {\em following} 
firing will cause, and $g(1, \epsilon)$ is the phase shift that an input immediately {\em preceding} firing will cause. 
However, we will now argue that in the neuroscience context, $g(0,\epsilon)$ and $g(1,\epsilon)$ should in fact be the same, and both should be zero.

The assumption
\begin{equation}
\label{eq:g0}
g(0,\epsilon)=0
\end{equation}
is justified by analogy with a neuron, which has a refractory period immediately following firing, during which
it is input-insensitive.
On the other hand, since input has little effect once the spike-generating currents have been activated, 
it is also reasonable to assume
\begin{equation}
\label{eq:g1}
g(1,\epsilon)=0.
\end{equation}
Finally, we assume that $g$ is smooth: 
\begin{equation}
\label{eq:smooth}
g \in C^1([0,1] \times [0,\infty)).
\end{equation}

We have defined phase as a number in $[0,1]$. It is mathematically tempting to view 
it instead as a point on the 
unit circle $S^1$, identifying $\varphi$ with the complex number $e^{2 \pi i \varphi}$. In that case, the mapping
$$
\varphi \mapsto \varphi + g(\varphi,\epsilon)
$$
(for fixed $\epsilon$)
becomes the mapping
$$
S^1 \rightarrow S^1
$$
described by 
$$
e^{2 \pi i \varphi} \mapsto e^{2 \pi i \varphi} e^{2 \pi i g(\varphi,\epsilon)}.
$$
Assumptions (\ref{eq:range_of_phi})--(\ref{eq:smooth}) make this mapping  continuous. However, continuity
of the mapping $S^1 \rightarrow S^1$ is a weaker
condition than (\ref{eq:range_of_phi})--(\ref{eq:smooth}); 
it would only require that the periodic extension of $g(\varphi,\epsilon)$ (fixed $\epsilon \geq 0$) from $[0,1)$ to $\mathbb{R}$ be continuous up to jumps of integer magnitude. Our more specific assumptions about $g$ are motivated by neuroscience, as I have explained.

\subsection{Two interacting pulse-coupled phase oscillators} 

\vskip 5pt
Consider two identical phase oscillators $A$ and $B$ with period $T$.
When $A$ reaches phase $1$ (fires), it sends a pulse to $B$, which causes the phase $\varphi_B$ of $B$
to shift from $\varphi_B$ to $\varphi_B + g(\varphi_B, \epsilon)$. Similarly, when $B$ reaches phase 1 (fires), it sends
a pulse to $A$, causing the shift $\varphi_A \mapsto \varphi_A + g(\varphi_A,\epsilon)$. 

We can now ask questions about the dynamics of $A$ and $B$. For instance we can ask whether synchrony 
is attracting, that is, whether the two 
oscillators will synchronize their firing if started out sufficiently close in phase to each other. 
}

Consider a time when $\varphi_A=0$ and $\varphi_B \in (0,1)$. So $A$ has just fired. We assume 
that it has had its effect on $B$ already. 
After time $(1 - \varphi_B)T$ elapses, $B$ fires. 
At that time, $A$ is at phase $1- \varphi_B$ and is advanced  to phase $1 - \varphi_B + g(1-\varphi_B,\epsilon)$. Now the 
roles reverse: $B$ is at phase 0 and $A$ is at a phase $1 - \varphi_B + g(1-\varphi_B,\epsilon) \in (0,1)$. After $A$
fires next and has its effect on $B$, $B$ will be at phase $F(\varphi_B,\epsilon)$ with 
\begin{equation}
\label{eq:defF}
F(\varphi,\epsilon) = 
 \varphi - g(1-\varphi,\epsilon) +
g(\varphi - g(1 - \varphi,\epsilon),\epsilon).
\end{equation}
We call $F$ the {\em strobe map} --- it is as though a strobe light were turned on each time $A$ completes
a cycle, to check where $B$ is. This idea for analyzing synchronization was suggested by Peskin in \cite{peskin_1975}.
I first saw the word ``strobe" for it in \cite{strogatz_sync}.

Since $F(\varphi,\epsilon)$ is a phase, we must always have $0 \leq F(\varphi,\epsilon) \leq 1$.
Indeed this follows easily from (\ref{eq:range_of_phi}). 
The two oscillators implement fixed point iteration for the function $F(\cdot,\epsilon)$:

{
\begin{equation}
\label{eq:fp_iteration}
\varphi_{k+1} = F(\varphi_k, \epsilon), ~~~ k=0,1,\ldots
\end{equation}
where $\varphi_k$ is the phase of $B$ after $A$ has completed its $k$-th cycle and has had
its effect on $B$. Therefore the analysis of synchronization now simply becomes an analysis of a fixed point iteration. 
Our assumptions imply that $\varphi = 0$ and $ \varphi = 1$ are fixed points: 
$$
F(0,\epsilon)=0, ~~~ F(1,\epsilon)=1.
$$
Synchrony is attracting if $0$ and $1$ are attracting fixed points. 

In the standard analysis of fixed point iteration, derivatives at 
fixed points plays a central role. As should be
expected, the derivatives of $F(\cdot, \epsilon)$ at $\varphi=0$ and $\varphi=1$ are
the same. A straightforward calculation shows the following lemma. 
\begin{lemma}  
\label{lemma:stable} 
Assume
(\ref{eq:nichts})--(\ref{eq:defF}). Then 
\begin{equation}
\label{eq:dF}
\frac{\partial F}{\partial \varphi}(0,\epsilon) = \frac{\partial F}{\partial \varphi}(1,\epsilon) = 
\left( 1 + \frac{\partial g}{\partial \varphi}(0,\epsilon) \right) \left( 1 + \frac{\partial g}{\partial \varphi}(1,\epsilon) \right).
\end{equation}
\end{lemma}

The idea of the strobe map is  usually attributed to Winfree \cite{Winfree_old}, although Peskin's 1975 course notes \cite{peskin_1975}, where the same idea appeared in a special context,
precede \cite{Winfree_old}.  Peskin's result
was later generalized by Mirollo and Strogatz in \cite{mirollo_strogatz_1990}. 
The literature on phase response functions and their 
use for the analysis of networks of oscillators has since become vast. For a review of phase response functions in neuroscience, 
see \cite{PRC_review}. Phase response functions have also been used in engineering; see for instance \cite{Kato_et_al_2019,Loe_et_al,Wang_et_al_2013}.

\subsection{Infinitesimal phase response function} 

\vskip 5pt

{ 
The linearization of $g$  around $\epsilon=0$, 
\begin{equation}
\label{eq:IPRC} 
\tilde{g}(\varphi,\epsilon) = \frac{\partial g}{\partial \epsilon}(\varphi,0) \epsilon,
\end{equation}
is commonly called the {\em infinitesimal phase response function}.

\begin{lemma} \label{lemma:is_PRC} Let $g \in C^2([0,1] \times [0,\infty))$ satisfy assumptions (\ref{eq:nichts})--(\ref{eq:g1}). 
Then for sufficiently small $\epsilon > 0$, 
(\ref{eq:nichts})--(\ref{eq:smooth}) hold with $g$ replaced by $\tilde{g}$.
\end{lemma}

\vskip 5pt
\noindent
{\bf \em Proof.} Because of the factor of $\epsilon$ in (\ref{eq:IPRC}), $\tilde{g}(\varphi,0) = 0$ for all $\varphi \in [0,1]$. Because $g \equiv 0$ for $\varphi=0$, 
also $\frac{\partial g}{\partial \epsilon} \equiv 0$ for $\varphi=0$, and
the analogous 
holds for $\varphi=1$. Therefore $\tilde{g}(0,\epsilon) = \tilde{g}(1,\epsilon)= 0$ for all $\epsilon \geq 0$. Since $g \in C^2([0,1] \times [0,\infty))$, $\tilde{g} \in C^1([0,1] \times [0,\infty))$. Finally, these properties imply that (\ref{eq:range_of_phi}) and (\ref{eq:range_of_phi_strict}) with $g$ replaced by $\tilde{g}$ hold for sufficiently small $\epsilon>0$. \qed

\vskip 5pt
The infinitesimal phase response function $\tilde{g}$ gives rise to the strobe map
$$
\tilde{F}(\varphi,\epsilon) = 
 \varphi - \tilde{g}(1-\varphi,\epsilon) +
\tilde{g}(\varphi - \tilde{g}(1 - \varphi,\epsilon),\epsilon).
$$

\section{Conditions for synchrony to be attracting} 

Throughout this section we consider two identical phase oscillators $A$ and $B$, interacting via a phase response
function $g \in C^2([0,1] \times [0,\infty))$ that satisfies (\ref{eq:nichts})--(\ref{eq:g1}). We assume that $F$ is the
strobe map defined by (\ref{eq:defF}). We take $\tilde{g}$ to be the associated infinitesimal phase response
function; it is a phase response function for sufficiently small $\epsilon>0$ by Lemma \ref{lemma:is_PRC}. 
The phrase ``for sufficiently small $\epsilon>0$" will appear often in the subsequent discussion; we will 
abbreviate it by ``for small $\epsilon$" from here on.

\subsection{Strongly and weakly attracting synchrony} 

\vskip 5pt
\begin{definition}
We say that synchrony is {\em $g$-attracting} if $0$ and $1$ are locally attracting
fixed points of $F(\cdot, \epsilon)$. We define  {\em $g$-repelling} analogously. 
\end{definition}

\vskip 5pt
The standard theory of fixed point iteration, together with eq.\ (\ref{eq:dF}), shows:

\begin{lemma}
Synchrony is $g$-attracting if 
\begin{equation}
\label{eq:strongly_attract}
\left| \left( 1 + \frac{\partial g}{\partial \varphi}(0,\epsilon) \right) \left( 1 + \frac{\partial g}{\partial \varphi}(1,\epsilon) \right) \right| < 1, 
\end{equation}
and it is $g$-repelling if 
\begin{equation}
\label{eq:strongly_repell}
\left| \left( 1 + \frac{\partial g}{\partial \varphi}(0,\epsilon) \right) \left( 1 + \frac{\partial g}{\partial \varphi}(1,\epsilon) \right) \right| > 1.
\end{equation}
\end{lemma}

It should be noted though that (\ref{eq:strongly_attract}) and (\ref{eq:strongly_repell}) are {\em sufficient} 
conditions for synchrony to be $g$-attracting and $g$-repelling, respectively, not {\em necessary} ones.

\begin{definition}
We say that synchrony is {\em strongly $g$-attracting} if (\ref{eq:strongly_attract}) holds. We say that
synchrony is {\em weakly $g$-attracting} if it is $g$-attracting but not strongly $g$-attracting. The terms {\em strongly $g$-repelling} and {\em weakly $g$-repelling} are defined analogously.
\end{definition}

To say that synchrony is strongly $g$-attracting is to say that the {\em linearized} analysis 
of the fixed point iteration shows that
$0$ and $1$ are attracting fixed points of $F$. To avoid confusion, we stress that this refers
to linearization of $F$ around $\varphi=0$ and $\varphi=1$; $\epsilon>0$ is fixed here. By contrast, $\tilde{g}$ is obtained by linearizing $g$ around $\epsilon=0$.

When synchrony is strongly $g$-attracting, two oscillators $A$ and $B$ interacting through the phase 
response function $g$ converge to synchrony exponentially fast when started out sufficiently close to synchrony.

Since $\tilde{g}$ is, for small $\epsilon$, a phase response function, we can also talk about 
synchrony being {\em $\tilde{g}$-attracting}, {\em strongly $\tilde{g}$-attracting}, and so on. 

\begin{lemma} \label{lemma:strongly_g_tilde_attracting} Synchrony is strongly $\tilde{g}$-attracting
for small $\epsilon$ if and only if one of the following two conditions 
holds. 
\begin{equation}
\label{one}
\frac{\partial^2 g}{\partial \varphi \partial \epsilon}(0,0) + \frac{\partial^2 g}{\partial \varphi \partial \epsilon}(1,0) <0,
\end{equation}
or 
\begin{equation}
\label{two}
\frac{\partial^2 g}{\partial \varphi \partial \epsilon}(0,0) = - \frac{\partial^2 g}{\partial \varphi \partial \epsilon}(1,0) \neq 0.
\end{equation}
Synchrony is strongly $\tilde{g}$-repelling for small $\epsilon$ if and only if
\begin{equation}
\label{three}
\frac{\partial^2 g}{\partial \varphi \partial \epsilon}(0,0) + \frac{\partial^2 g}{\partial \varphi \partial \epsilon}(1,0) >0,
\end{equation}
\end{lemma}

\vskip 5pt
\noindent
{\bf \em Proof.} By definition, synchrony is strongly $\tilde{g}$-attracting if and only if 
$$
\left| \left( 1 + \frac{\partial \tilde{g}}{\partial \varphi}(0,\epsilon) \right) \left( 1 + \frac{\partial \tilde{g}}{\partial \varphi}(1,\epsilon) \right) \right| < 1. 
$$
By definition of $\tilde{g}$, this means 
$$
\left| \left( 1 + \frac{\partial^2 g}{\partial \varphi \partial \epsilon}(0,0)  \epsilon \right) 
 \left( 1 + \frac{\partial^2 g}{\partial \varphi \partial \epsilon}(1,0)  \epsilon \right)  \right| < 1, 
$$
that is, 
$$
\left| 1 + \left( \frac{\partial^2 g}{\partial \varphi \partial \epsilon}(0,0)   + \frac{\partial^2 g}{\partial \varphi \partial \epsilon}(1,0)   \right) \epsilon + \frac{\partial^2 g}{\partial \varphi \partial \epsilon}(0,0)   
\frac{\partial^2 g}{\partial \varphi \partial \epsilon}(1,0)   \epsilon^2 \right| < 1.
$$
This clearly holds for sufficiently small $\epsilon>0$ if and only if either (\ref{one}) or (\ref{two}) holds. 

Similarly, synchrony is strongly $\tilde{g}$-repelling if and only if 
$$
\left| 1 + \left( \frac{\partial^2 g}{\partial \varphi \partial \epsilon}(0,0)   + \frac{\partial^2 g}{\partial \varphi \partial \epsilon}(1,0)   \right) \epsilon + \frac{\partial^2 g}{\partial \varphi \partial \epsilon}(0,0)   
\frac{\partial^2 g}{\partial \varphi \partial \epsilon}(1,0)   \epsilon^2 \right| > 1.
$$
This holds for sufficiently small $\epsilon>0$ if and only if (\ref{three}) holds. 
\qed

\subsection{Why weakly attracting synchrony may not be unusual} 

\vskip 5pt
 If inputs arriving around the time of firing have very little effect --- as
is typical for a neuron --- then it seems not unreasonable that (\ref{eq:g0}) and (\ref{eq:g1}) could
be strengthened like this:
\begin{equation}
\label{eq:gg0}
g(0,\epsilon)= \frac{\partial g}{\partial \varphi}(0,\epsilon) = 0,
\end{equation}
and
\begin{equation}
\label{eq:gg1}
g(1,\epsilon)= \frac{\partial g}{\partial \varphi}(1,\epsilon) = 0.
\end{equation}
For example (\ref{eq:gg0}) clearly holds if there is an {\em absolute} refractory period, a brief time following firing during
which the neuron is {\em entirely} input-insensitive. 
Both (\ref{eq:gg0}) and (\ref{eq:gg1}) hold for an Ermentrout-Kopell theta neuron responding to
 instantaneous charge injections; see Section 4. 
 Note that (\ref{eq:gg0}) and (\ref{eq:gg1}) imply that neither (\ref{eq:strongly_attract}) nor (\ref{eq:strongly_repell}) hold.

\subsection{Very strongly $\tilde{g}$-attracting synchrony} 

\vskip 5pt
\begin{definition}
We say that synchrony is {\em very strongly $\tilde{g}$-attracting for small $\epsilon$} if (\ref{one}) holds. 
\end{definition}

We don't define {\em very strongly $\tilde{g}$-repelling}, nor do we define {\em very strongly $g$-attracting} or {\em 
very strongly $g$-repelling}.

\section{Results} 

\begin{theorem}
It is possible for synchrony to be weakly $\tilde{g}$-attracting for small $\epsilon$, yet not $g$-attracting
for any $\epsilon$. Two identical Ermentrout-Kopell theta neurons interacting via instantaneous charge
injection provide an example. (See Section 4 for details.) 
\end{theorem}

Theorem 2 strengthens Theorem 1.

\begin{theorem}
It is possible for synchrony to be weakly, or even strongly $\tilde{g}$-attracting, yet strongly $g$-repelling for
all $\epsilon$.
\end{theorem}

Theorem 3 provides conditions under which the infinitesimal phase response function does make
the correct prediction.

\begin{theorem} (a) 
If synchrony is very strongly $\tilde{g}$-attracting for small $\epsilon$,  it is strongly $g$-attracting 
for small $\epsilon$. (b) If synchrony is strongly $\tilde{g}$-repelling for small $\epsilon$, it is strongly $g$-repelling for small $\epsilon$.
\end{theorem}

\section{Proof of Theorem 1}

\subsection{The theta neuron} 
\vskip 5pt

Ermentrout and Kopell \cite{ermentroutkopell} proposed to model a neuron as a point moving along the unit circle, 
with its position $(\cos \theta(t), \sin \theta(t))$ governed by 
\begin{equation}
\label{eq:theta_neuron}
\frac{d \theta}{dt} = 1 - \cos \theta + I (1 + \cos \theta).
\end{equation}
The right-hand side  of (\ref{eq:theta_neuron}) is positive for all $\theta$ if $I>0$, 
which we will assume here. The period, that is, the time it takes for $\theta$ to increase by $2 \pi$ and therefore for the moving
point to move once around the unit circle, is then 
\begin{equation}
\label{eq:period_theta_neuron}
T = \int_{-\pi}^\pi \frac{dt}{d \theta} ~\! d \theta =  \int_{-\pi}^\pi \frac{1}{1-\cos \theta + I (1+ \cos \theta) } ~\! d \theta = \frac{\pi}{\sqrt{I}}.
\end{equation}
See  \cite{Ermentrouttheta}, \cite{borgerskopell2003}, or \cite{borgers_book} for discussions of what (\ref{eq:theta_neuron}) has to do with a neuron, and also for
an explanation why 
\begin{equation}
\label{eq:v_of_theta}
v = \frac{1}{2} + \frac{1}{2} \tan \frac{\theta}{2}
\end{equation}
should be considered the analogue of the ``membrane potential" (interior voltage) of the neuron --- notwithstanding the fact that it becomes infinite
when $\theta$ is an odd multiple of $\pi$. We say that the theta neuron ``fires" when $\theta$ reaches an odd multiple of $\pi$. 

\subsection{Instantaneous charge injections} 
\vskip 5pt

An instantaneous injection of a positive amount of charge into a neuron would make the membrane potential jump. We will 
therefore assume that the response of a theta neuron to a brief input pulse is the jump 
\begin{equation}
\label{eq:v_jump}
v \mapsto v + \Delta v.
\end{equation}
We'll assume $\Delta v>0$ here ({\em positive} charge is injected, the input is {\em excitatory}). 
Using eq.\ (\ref{eq:v_of_theta}), eq.\ (\ref{eq:v_jump}) translates into 
\begin{equation}
\label{eq:theta_jump}
\theta \mapsto 
2 \arctan \left(   \tan \frac{\theta}{2} + 2\Delta v \right).
\end{equation}

\subsection{The phase response function} 
\vskip 5pt

Equation (\ref{eq:theta_jump}) 
does not  yet describe the {\em phase} response, since $\theta$ is not the same as the phase $\varphi$.
The relation between $\theta$ and $\varphi$ is
\begin{equation}
\label{eq:theta_to_phi}
\varphi = 
\frac{1 - \int_{\theta}^\pi \frac{1}{1-\cos s + I(1+\cos s)} d s }{T}  = \frac{1}{2} + \frac{1}{\pi} \arctan \left( \frac{\tan \frac{\theta}{2}}{\sqrt{I}} \right).
\end{equation}
We use (\ref{eq:theta_to_phi}) to express the right-hand side of (\ref{eq:theta_jump}) in terms of $\varphi$, and conclude { that (\ref{eq:theta_jump}) translates into}
$$
\varphi \mapsto \varphi +  \frac{1}{\pi} \arctan \left( \tan \left( \left( \varphi- \frac{1}{2} \right) \pi \right) + \frac{2 \Delta v}{\sqrt{I}} \right) - 
\left( \varphi - \frac{1}{2} \right).
$$
To simplify the notation, we write 
$$
\epsilon = \frac{2 \Delta v}{\sqrt{I}}.
$$
So 
\begin{equation}
\label{eq:def_g}
g(\varphi,{ \epsilon}) = \frac{1}{\pi} \arctan \left( \tan \left( \left( \varphi- \frac{1}{2} \right) \pi \right) +{ \epsilon} \right)  - 
\left( \varphi - \frac{1}{2} \right).
\end{equation}

\noindent
A straightforward computation verifies  (\ref{eq:gg0}) and (\ref{eq:gg1}) now. 

\subsection{Synchrony is neutrally $g$-stable for all $\epsilon \geq 0$}

\vskip 5pt
We insert (\ref{eq:def_g})  into  (\ref{eq:defF}) and find a surprise: 
\begin{equation}
\label{eq:F_is_the_identity}
F(\varphi, { \epsilon}) = \varphi ~~~\mbox{for all $\varphi \in [0,1]$, ~$\epsilon \geq 0$}.
\end{equation}
The exchange of one pulse from $B$ to $A$ and one from $A$ to $B$ brings the phase difference between the two theta neurons back to where it started. Any initial phase difference persists. This is true for all $\epsilon>0$.

\subsection{Synchrony is weakly $\tilde{g}$-attracting for small $\epsilon$} 

\vskip 5pt
We now replace (\ref{eq:def_g}) by its local linear approximation near { $\epsilon=0$:}
\begin{equation}
\label{eq:infinitesimal_PRC}
\tilde{g}(\varphi,{ \epsilon}) = \frac{1}{\pi} ~\!  \frac{1}{1 + \tan^2 \left( \left( \varphi - \frac{1}{2} \right) \pi \right) }  ~\! \epsilon.
\end{equation}
For sufficiently small $\epsilon>0$, for example for $\epsilon \leq 1$,  $\tilde{g}<1-\varphi$ for $\varphi \in (0,1)$, so $\tilde{g}$ is indeed a valid
phase response function. 
We insert (\ref{eq:infinitesimal_PRC}) into (\ref{eq:defF}) to find the strobe map
\begin{equation}
\label{eq:infinitesimal_F}
\tilde{F}(\varphi,{ \epsilon}) = \varphi+ \frac{\epsilon}{2 \pi} \cos(2 \pi \varphi) - \frac{\epsilon}{2 \pi } \cos \left( 2 \pi \varphi + \epsilon \cos(2 \pi \varphi) - \epsilon \right).
\end{equation}

Let $\epsilon>0$ be fixed. 
A straightforward calculation shows that in the limit as $\varphi \rightarrow 0$, 
\begin{equation}
\label{eq:tilde_F_expansion}
\tilde{F}(\varphi,\epsilon) = 
\varphi - 2 \epsilon^2 \pi^2 \varphi^3 + O(\varphi^4) < \varphi.
\end{equation}
This implies that for sufficiently small $\varphi$, $\tilde{F}(\varphi,\epsilon)< \varphi$. Therefore the sequence
$\varphi_0, \varphi_1,\ldots$ generated by fixed point iteration is decreasing if $\varphi_0$ is
close enough to $0$, and since it is bounded below by $0$
it must converge. The limit must be a fixed point of $\tilde{F}(\cdot, \epsilon)$, and for sufficiently small 
$\varphi_0>0$, that fixed point must be $0$. 
So 
 $\varphi=0$ is an attracting fixed point.

Arguments like those given in the preceding paragraph show that 
$\varphi=1$ is an attracting fixed point as well. 
We conclude that 
synchrony is $\tilde{g}$-attracting. It is {\em weakly} $\tilde{g}$-attracting because 
$\frac{\partial \tilde{F}}{\partial \varphi}(0,\epsilon) = \frac{\partial \tilde{F}}{\partial \varphi}(1,\epsilon) =1$. 
On the other hand, when the phase response function is $g$, synchrony is just {\em neutrally} stable, not 
attracting. 

This proves Theorem 1.

\section{Proof of Theorem 2} 

\vskip 5pt
\subsection{Example 1: Weakly $\tilde{g}$-attracting, yet strongly $g$-repelling}

\vskip 5pt
This example is a modification of the phase response function of the theta
neuron:

\begin{equation}
\label{eq:add_this}
g(\varphi,\epsilon) = \frac{1}{\pi} \arctan \left( \tan \left( \left( \varphi- \frac{1}{2} \right) \pi \right) +{ \epsilon} \right)  - 
\left( \varphi - \frac{1}{2} \right) + \varphi (1-\varphi)^{ 2} \epsilon^2.
\end{equation}
Compare this with eq.\ (\ref{eq:def_g}): It is the phase response function for the theta neuron, plus $\varphi(1-\varphi)^2 \epsilon^2$.
Unfortunately 
I am not aware of any neuroscience motivation for this formula.
I added the term 
$\varphi(1-\varphi)^{ 2} \epsilon^2$ simply to construct an interesting mathematical example. 

Using that (\ref{eq:gg0}) and (\ref{eq:gg1}) hold for the phase response function of the theta neuron, we have 

$$
\left( 1 + \frac{\partial {g}}{\partial \varphi}(0,\epsilon) \right) \left( 1 + \frac{\partial {g}}{\partial \varphi}(1,\epsilon) \right) = 1+\epsilon^2>1.
$$
So synchrony is strongly ${g}$-repelling for all $\epsilon>0$.
However, the {\em infinitesimal} phase response function, namely the linearization 
 around $\epsilon=0$, is the same as for the theta neuron --- we just added a term of size $O(\epsilon^2)$. Synchrony is weakly $\tilde{g}$-attracting for small $\epsilon$, but strongly $g$-repelling for all $\epsilon>0$.
 
 \subsection{Example 2: Strongly $\tilde{g}$-attracting yet strongly $g$-repelling} 
 \vskip 5pt
 
Let 
\begin{equation}
\label{g_3}
g(\varphi,\epsilon) = \varphi(1-\varphi) \epsilon - 2 \varphi(\varphi-1)^2 (2 \varphi-1) \epsilon^2.
\end{equation}
Here 
$$
\frac{\partial^2 g}{\partial \varphi \partial \epsilon}(0,0) =  - \frac{\partial^2 g}{\partial \varphi \partial \epsilon} (1,0) = 1.
$$
Therefore synchrony is strongly $\tilde{g}$-attracting for small $\epsilon$ by Lemma \ref{lemma:strongly_g_tilde_attracting}. However, 
$$
\left( 1 + \frac{\partial g}{\partial \varphi}(0,\epsilon) \right) \left( 1 + \frac{\partial g}{\partial \varphi}(1,\epsilon) \right)  
=  \left( 1 + \epsilon + 2 \epsilon^2 \right) \left( 1 - \epsilon \right)  = 1 +  \epsilon^2 -2 \epsilon^3>1
$$
for small $\epsilon$, so synchrony is strongly $g$-repelling for small $\epsilon$.

\vskip 5pt
This concludes the proof of Theorem 2.

\section{Proof of Theorem 3} 

Now suppose that synchrony is {\em very} strongly $\tilde{g}$-attracting. Recall that by definition, this means
that (\ref{one}) holds. Then, expanding around $\epsilon=0$, we find
\begin{eqnarray*}
&~&
\left( 1 + \frac{\partial g}{\partial \varphi}(0,\epsilon) \right) \left( 1 + \frac{\partial g}{\partial \varphi}(1,\epsilon) \right)   \\ 
&=&
\left( 1 + \frac{\partial^2 g}{\partial \varphi \partial \epsilon}(0,0) \epsilon + o(\epsilon) \right) 
\left( 1 + \frac{\partial^2 g}{\partial \varphi \partial \epsilon}(1,0) \epsilon + o(\epsilon) \right)  \\
&=&1 +  \left( \frac{\partial^2 g}{\partial \varphi \partial \epsilon}(0,0)  + \frac{\partial^2 g}{\partial \varphi \partial \epsilon}(1,0)  \right) \epsilon + o(\epsilon) < 1
\end{eqnarray*}
for small $\epsilon$. This proves  part (a) of Theorem 3. Part (b) is proved analogously.

{
\section{Summary and a concluding comment} 

\subsection{Summary} 
\vskip 5pt
The infinitesimal phase response function first approximates the behavior of $g$ near $\epsilon=0$, then 
analyzes synchrony by studying the behavior near $\varphi=0$ and $\varphi=1$. What one really wants to know is the behavior near $\varphi=0$ and $\varphi=1$ for fixed small $\epsilon>0$. There is no reason why these should
be the same, and the examples given here show that they aren't always the same. 

\subsection{Comparison with the Hartman-Grobman theorem} 

\vskip 5pt

The Hartman-Grobman theorem states that if linearization near an equilibrium of a dynamical system yields a
definite prediction about the dynamics near the equilibrium, that prediction holds (qualitatively) for the non-linear problem. 
See \cite{Guysinsky:2003p20294} for a particularly nice version of the theorem.
The examples I have given
here show that if linearization around $\epsilon=0$  yields
a definite prediction about whether synchrony is attracting, that prediction may still be false 
for the non-linear problem, for arbitrarily small but possible interaction strengths. The analogue of the Hartman-Grobman theorem that one might expect at first sight 
does not hold here.

Theorem 2 shows that even if linearization of $g$ around $\epsilon=0$, then linearization of $F$ around 
$\varphi=0$ and $\varphi=1$, yields a definite prediction for small $\epsilon>0$, that prediction may still be false.
In other words, the analogue of the Hartman-Grobman theorem that one might expect at {\em second} sight 
is also false.

Theorem 3 gives a condition under which the prediction of the analysis based on $\tilde{F}$ {\em does} hold for
$F$, but Section 2.2 explains why situations in which that condition is satisfied may not be typical.

}

\vskip 10pt
\noindent
{\bf Acknowledgments.} I thank an anonymous reviewer for an exceptionally thorough and thoughtful report, which led to very substantial clarifications, improvements, and extensions.  I also thank Nancy Kopell, whose 
comments and questions
resulted in Theorems 2 and 3.





%
%
%

\bibliographystyle{elsarticle-num}
 
\end{document}